\documentclass{amsart}
\usepackage{graphicx}
\usepackage{amssymb,amscd,amsthm,amsxtra}
\usepackage{hyperref}
\usepackage{latexsym}
\usepackage{epsfig}
\usepackage{mathtools}
\usepackage{esint}
\usepackage{color}

\newtheorem{thm}{Theorem}[section]
\newtheorem{cor}[thm]{Corollary}

\theoremstyle{definition}
\newtheorem{defn}[thm]{Definition}
\theoremstyle{remark}
\newtheorem{rem}[thm]{Remark}
\numberwithin{equation}{section}

\newcommand{\R}{\mathbb R}

\newcommand{\eps}{\epsilon}

\newcommand{\p}{\partial}

\newcommand{\comment}[1]{}

\begin{document}

\title[Boundary Harnack Inequality]{A short proof of Boundary Harnack Inequality}
\author{D. De Silva}
\address{Department of Mathematics, Barnard College, Columbia University, New York, NY 10027}
\email{\tt  desilva@math.columbia.edu}
\author{O. Savin}
\address{Department of Mathematics, Columbia University, New York, NY 10027}\email{\tt  savin@math.columbia.edu}
\begin{abstract}We give a direct analytic proof of the classical Boundary Harnack inequality for solutions to linear uniformly elliptic equations in either divergence or non-divergence form.
 \end{abstract}

\maketitle
\section{Introduction}

In this note we give a short proof of the classical Boundary Harnack inequality for solutions to linear uniformly elliptic equations, which is based only on the Harnack inequality and Weak Harnack inequality. Our proof provides a unified approach for both divergence and non-divergence linear equations. The strategy applies to several other extensions of the Boundary Harnack inequality which we mention in the last section. 

We recall the setting of the classical Boundary Harnack inequality in Lipschitz domains. 
As usual, we write $x=(x',x_n) \in \R^n$ and $B'_1$ denotes the unit ball in $\R^{n-1},$ centered at $0.$

Given a Lipschitz function $g$, with $$g: B'_1 \subset \R^{n-1} \to \R, \quad g \in C^{0,1}, \quad \|g\|_{C^{0,1}} \leq L, \quad  g(0)=0,$$ 
we denote by
$$\Gamma:= \{x_n =g(x')\}, \quad h_\Gamma(x):= x_n -g(x').$$
We define the cylindrical region of radius $r$ and height $\rho$ above the graph $\Gamma$ as
$$\mathcal C(r,\rho):=\{(x', x_n) \ : \ |x'| <r, 0 < h_\Gamma <  \rho\};$$ if $r=\rho$ we write simply $\mathcal C_r:=\mathcal C(r,r).$ 
Let $\mathcal Lu$ be either
$$\mathcal L u : = div(A(x) \nabla u), \quad \text{or} \quad \mathcal L u := tr(A(x)D^2 u)$$ with $A$ uniformly elliptic, that is 
$$\lambda I \leq A \leq \Lambda I,\quad 0<\lambda \leq \Lambda <+\infty.$$
Boundary Harnack Inequality states the following.
\begin{thm}\label{main} Let $u,v>0$ satisfy $\mathcal L u=\mathcal L v=0$ in $\mathcal C_1$ and vanish continuously on $\Gamma$. Assume $u,v$ are normalized so that $u\left( e_n/2\right)=v(e_n/2)=1,$ then
\begin{equation}\label{BHI} C^{-1} \leq \frac u v \leq C, \quad \text{in $\mathcal C_{1/2},$}\end{equation} with $C$ depending on $n, \lambda, \Lambda,$
and $L$. \end{thm}

The classical case when $\mathcal L= \Delta$ appears in \cite{A,D,K,W}. Operators in divergence form were first considered in \cite{CFMS} while the case of operator in non-divergence form was treated in \cite{FGMS}. The same result for operators in divergence form was extended also to so-called NTA domains in \cite{JK}. The case of H\"older domains and $\mathcal L$ in divergence form was  addressed with probabilistic techniques in \cite{BB1,BBB}, and an analytic proof was then provided in \cite{F}. For H\"older domains and operators $\mathcal L$ in non-divergence form, it is necessary that the domain is $C^{0,\alpha}$ with $\alpha>1/2$ (or that it satisfies a uniform density property), and the proof again is based on a probabilistic approach \cite{BB2}. 

The purpose of this short note is to provide a unified analytic proof of Theorem \ref{main} that does not make use of the Green's function and which holds for both operators in non-divergence and in divergence form. The idea is to find an ``almost positivity property" of a solution, which can be iterated from scale 1 to all smaller scales. The same strategy also applies to other similar situations like that of H\"older domains, NTA domains and to the case of ``slit" domains (see Section 3, for the precise definition of NTA and slit domain.) 

The key property of uniformly elliptic equations needed in our proof is the following Weak Harnack Inequality for subsolutions, which holds for both divergence  \cite{DG} and non-divergence \cite{KS} equations.

\begin{thm}\label{WH} Let $\mathcal L \varphi \geq 0$ in $B_1,$  with $\varphi \leq 1$ in $B_1.$ If, for some $\eta>0,$
$$|\{\varphi \leq 0\}| \geq \eta$$
then
$$ \varphi \leq  1-c(\eta) \quad \quad \mbox{in $B_{1/2}$},$$
with $0<c(\eta)<1$ depending on $\eta$ and $n,\lambda, \Lambda.$
\end{thm}

A consequence of Theorem \ref{WH} based on scaling and covering arguments implies the following version of the Weak Harnack inequality (see for example Lemma 4.4 in \cite{CC} and Theorem 8.17 in \cite{GT}.)

\begin{thm}\label{lp}  Let $\mathcal L \varphi \geq 0$ in $B_1,$ then for every $p>0$
$$\|\varphi^+\|_{L^\infty(B_{1/2})} \leq C(p) \|\varphi\|_{L^p(B_1)},$$
with $C(p)>0$ depending on $p$ and $n, \lambda, \Lambda$.
\end{thm}

The paper is organized as follows. In Section 2 we provide the proof of Theorem \ref{main}. In Section 3 we provide several extensions, precisely  we discuss the case of H\"older domains, NTA domains, and slit domains.

\section{Proof of Theorem \ref{main}}



We divide the proof of Theorem \ref{main} in three steps. In what follows, constants depending only on $n, L,\lambda, \Lambda$ are called universal. Recall that,
$$\mathcal C(r,\rho):=\{(x', x_n) \ : \ |x'| <r, 0 < h_\Gamma < \rho\},$$ and if $r=\rho$ we write simply $\mathcal C_r.$ We also set, for $\delta>0$ small to be made precise later, $$\mathcal A_{r}:= \mathcal C_r \setminus \mathcal C(r,r\delta).$$
The idea is to show that a solution which is large in $\mathcal A_r$ and not too negative in $\mathcal C_r$ will remain positive in all smaller subdomains $\mathcal A_{2^{-k}r}$, with $k \ge 1$.

\

\textit{Step 1.} There exist $M>0$ large and $\delta>0$ small universal, such that if $w$ is a solution to $\mathcal L w=0$ in $\mathcal C_1$ (possibly changing sign) which vanishes on $\Gamma$, and
$$w \geq M \quad \text{on $\mathcal A_{1}$},$$ and
$$w \geq -1 \quad \text{on $\mathcal C_1$},$$
then,
\begin{equation}\label{201}
w \geq Ma \quad \text{on $\mathcal A_{\frac 1 2}$},
\end{equation} and
\begin{equation}\label{202}
w \geq -a \quad \text{on $\mathcal C_{\frac 1 2}$}, 
\end{equation}
for some small $a>0$. 

Notice that after a dilation we can apply the conclusion again to $w/a$. By rescaling and iterating, we conclude that $$w \ge M a_k \quad \quad \text{on $\mathcal A_{2^{-k}}$},$$
for a sequence of positive numbers $a_k$, hence
$$w > 0 \quad \text{on the line segment $\{te_n, \quad 0<t<1\}.$}$$

To prove our claim we first establish a lower bound for $w$ on $\mathcal A_{\frac 1 2}$. Pick $x_0$ in this set. The Lipschitz continuity of $g$ implies that the cone of vertex at $(x_0', g(x_0'))$ and slope $L$ is included in $\mathcal C_1$ near its vertex. We apply the interior Harnack inequality to $w+1 \ge 0$ in a sequence of overlapping balls included in this cone which connect $x_0$ with $x_0 + \delta e_n  \in \mathcal A_{1}$, and notice that the number of balls needed depends only on $L$. We conclude that 
\begin{equation}\label{hi}w(x_0) \geq (M+1)c_L -1,\end{equation}
for some $0<c_L<1$ universal. We choose $$a:=c_L/2,$$ and then $M$ large, universal, such that the right-hand-side in \eqref{hi} is larger than $Ma$. Thus \eqref{201} is established.

Next we show that \eqref{202} holds.
Let $x_0 \in \mathcal C(1-\delta, \delta)$, and let $Q_\delta(x_0)$ be a cylinder of height $2\delta$ and radius $\delta$ around $x_0$. Assume $w$ is extended to zero in $\{h_\Gamma<0\}.$ Notice that on each vertical segment of $Q_\delta(x_0)$, $w^-=0$ on at least half of its length (the parts which fall in $\mathcal A_{1}$ or below the graph of $g$). By Weak Harnack inequality applied to $w^-$ in $Q_\delta(x_0)$, and recalling that by hypothesis $w^- \le 1$ in $\mathcal C_1$, we conclude that $$w^-(x_0) \leq 1-c_0,$$ for some small $0<c_0<1$ universal. Hence $w^- \le 1-c_0$ in $\mathcal C_{1-\delta}$. By iterating this result we find $w^- \le (1-c_0)^k$ on $\mathcal C_{1-k \delta}$, hence
\begin{equation}\label{lip1}w \geq -(1-c_0)^\frac{1}{4\delta} \quad \text{on $\mathcal C_{\frac 1 2}.$}\end{equation}
Now we can choose $\delta>0$ small, universal, so that 
$$(1-c_0)^\frac{1}{4\delta} \leq a,$$
and Step 1 is proved.

\begin{figure}[h]
\includegraphics[width=0.6 \textwidth]{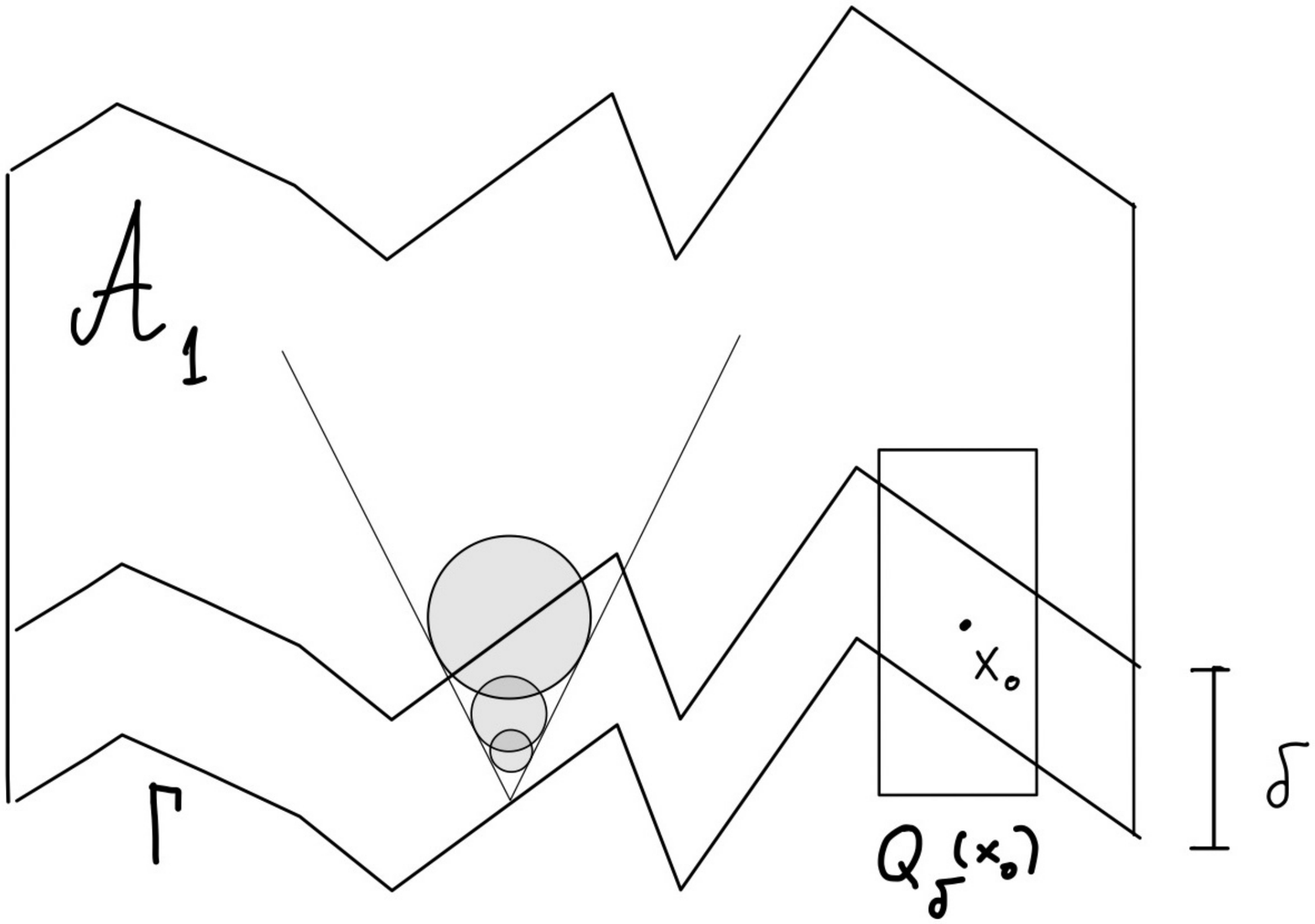}
\caption{Step 1}
    \label{fig1}
\end{figure}

\

\textit{Step 2.} In this step we show that 
\begin{equation} u,v \leq C \quad \text{in $\mathcal C_{1/2}$},\end{equation}
with $C>0$ universal.

We extend $u=0$ in $\{h_\Gamma<0\},$ and we still denote the resulting function by $u$. Then $\mathcal L u \geq 0$ say in $B'_1 \times \{-1 <h_\Gamma < 1\}$ and by Theorem \ref{lp}, for any $p>0$,
\begin{equation}\label{linf} \|u\|_{L^\infty(\mathcal C_{1/2})} \leq C(p,L)  \, \, \|u\|_{L^p(\mathcal C_{3/4})}.\end{equation}

On the other hand, since $g$ is Lipschitz, an iterated application of the interior Harnack inequality gives that
\begin{equation}\label{vert}
u(x) \leq h_\Gamma(x)^{-K}, \quad x \in \mathcal C_{3/4},
\end{equation}
for some large universal $K.$ Thus, $u^p$ is integrable by choosing $p=1/2K$, and our conclusion follows from \eqref{linf}. 

\





\textit{Step 3.} 
We show that for a large constant $C^*>0$ universal, 
$$w:=C^* u - v \geq 0 \quad \text{in $\mathcal C_{1/2}.$}$$
By Step 2 we know that $v \leq C_0$ in $\mathcal C_{3/4}.$ Moreover,
given $\delta>0$ universal from Step 1, since $u(e_n/2)=1$, we conclude by interior Harnack that
$$u \geq c_L(\delta), \quad \text{in $\mathcal A_{5/8}$}.$$ Thus, we can choose $C^*$ large so that
$$w \geq M C_0\quad \text{on $\mathcal A_{5/8}$},$$ with $M$ the constant in Step 1.
Moreover, 
$$w \geq -C_0 \quad \text{on $\mathcal C_{5/8}.$}$$
We conclude by Step 1 that $w \geq 0$ on the line $\{t e_n, 0<t<5/8\}$. By repeating the same argument at all points on $\Gamma \cap \overline {\mathcal C}_{1/2}$, our claim follows.

\qed

\section{Further extensions}

In this section we extend Theorem \ref{main} to other similar situations. The main point is to choose the domains $\mathcal C_r$ and $\mathcal A_r$ accordingly in each case, so that a quantified positivity statement as in Step 1 in Section 2 can be iterated.

\subsection{The case of H\"older domains.} Assume that $g \in C^{0,\alpha}$ with $\alpha>\frac 1 2.$ We prove here that the statement of Theorem \ref{main} remains valid. We mention that our proof shows that the constant $C$ does not depend on $\|v\|_{L^\infty}$, which is assumed in \cite{BB2} for the case of operators in non-divergence form.

The proof follows the same steps as in the Lipschitz case. We sketch below only Steps 1 and 2, as Step 3 is basically unchanged. Constants depending on $n,\alpha, \lambda, \Lambda$ and  $\|g\|_{C^{0,\alpha}}$ are now called universal. Here,
$$\mathcal A_r:=\mathcal C_r \setminus \mathcal C(r,r^\beta),$$
for $\beta>1$ to be made precise later.

\


\textit{Step 1.}  We show that,
there exist $C_0, \beta>1$ universal, such that if $w$ is a solution to $\mathcal L w=0$ in $\mathcal C_r$ (possibly changing sign) which vanishes on $\Gamma$ and
$$w \geq f(r) \quad \text{on $\mathcal A_{r}$}, $$
and
$$w \geq -1 \quad \text{on $\mathcal C_r$},$$
where $$f(r):= e^{C_0r^\gamma}, \quad \gamma:=\beta(1-\frac{1}{\alpha}) <0,$$
then,
\begin{equation}\label{301}
w \geq f( \frac r 2) \, \, a \quad \text{on $\mathcal A_{\frac r 2}$},
\end{equation} and
\begin{equation}\label{302}
w \geq -a \quad \text{on $\mathcal C_{\frac r 2}$}, 
\end{equation}
for some small $a=a(r)>0$, as long as $r\leq r_0$ universal.

The conclusion can be iterated and we obtain as before that if the hypotheses are satisfied in $\mathcal C_{r_0}$ then 
$$w > 0 \quad \text{on the line segment $\{te_n, \quad 0<t<r_0\}.$}$$

We argue as in the case of Lipschitz domains, however, since $g$ is H\"older continuous, when applying interior Harnack inequality to $w+1$, we need $$C(r^\beta)^{1-\frac 1 \alpha} =C r^\gamma \quad \mbox{balls}$$ to connect a point in $\mathcal A_{r/2}$ with a point in $\mathcal A_r$. We conclude that
\begin{equation}\label{w1}w \geq (f(r)+1) e^{-C_1r^\gamma} -1 \quad \text{in $\mathcal A_{r/2}$},\end{equation}
for some $C_1$ universal. We obtain for $r$ small, 
\begin{equation}\label{w1}w \ge 1 \quad \text{in $\mathcal A_{r/2}$},\end{equation} by choosing $C_0=2 C_1$.
 Moreover, as in the Lipschitz case (see \eqref{lip1} where we only used the graph property of $\Gamma$), 
$$w \geq - e^{-c_0 r^{1-\beta}} \quad \text{in $\mathcal C_{r/2}$},$$
for $c_0$ small universal.
We choose 
$$a(r):=e^{-c_0 r^{1-\beta}},$$
hence in view of \eqref{w1}, our claim 
$$w \ge 1 \ge a(r) f(\frac r2) $$is satisfied for $r$ small, as long as we can pick $\beta$ such that 
$$\gamma :=\beta(1-\frac 1 \alpha) > 1-\beta.$$ This is possible for $\alpha >1/2.$

\begin{figure}[h]
\includegraphics[width=0.6 \textwidth]{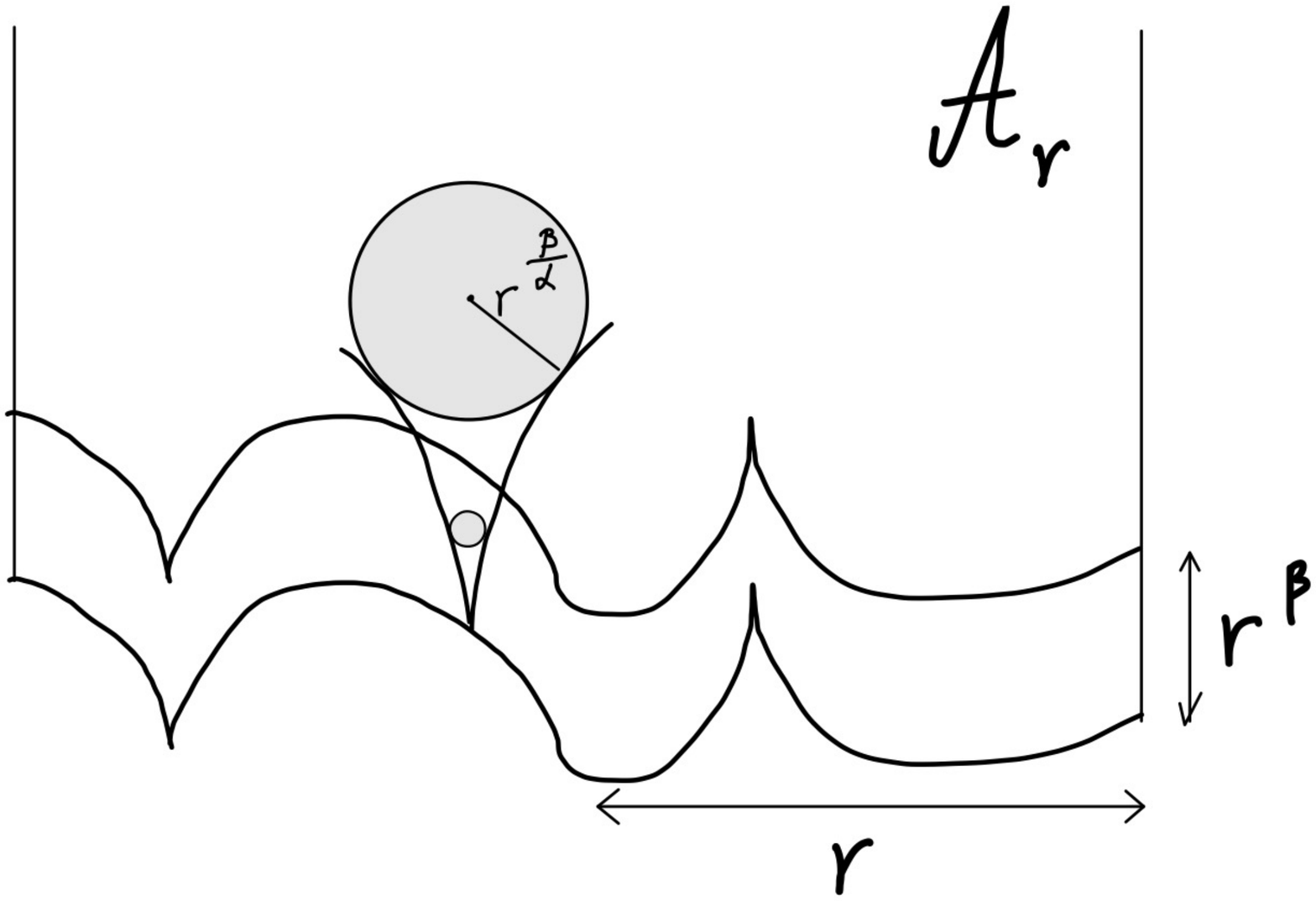}
\caption{Step 1}
    \label{fig2}
\end{figure}

\

\textit{Step 2.} We show that,
$$u,v \leq C_2 \quad \text{in $\mathcal C_{1/2}$,}$$
with $C_2$ universal. Here we apply an iterative argument similar to the one in Step 1 above, as \eqref{vert} no longer holds and our claim cannot be deduced by a direct application of Theorem \ref{lp}.
Since $u(e_n/2)=1$, the interior Harnack inequality gives that
\begin{equation}\label{int}u \leq e^{C_1 h_\Gamma^{1-1/\alpha}} \quad \quad \text{in $\mathcal C_{3/4}$},\end{equation}
 with $C_1$ universal.
With the same notation as Step 1, we wish to prove that if $r$ is smaller than a universal $r_0$ and
$$u(y) \geq f(r),$$ for some $y \in \mathcal C_{1/2}$,
then we can find $$z \in S:=\{|y'-z'|=r, \quad 0<h_{\Gamma}(z)< r^\beta\},$$ such that
$$u(z) \geq f(\frac{r}{2}).$$
Since $|z-y|\le C r^\beta$, we see that for $r$ small enough, we can build a convergent sequence of points $y_k \in \mathcal C_{3/4}$ with $u(y_k)\ge f(2^{-k}r )\to \infty$, a contradiction.

To show the existence of the point $z$ we let 
$$w:=(u-e^{C_1 r^\gamma})^+.$$ Then, in view of \eqref{int},
$$w = 0 \quad \text{on $T:=\{|x'-y'|\leq r, h_{\Gamma}= r^\beta\}$},$$
and $w=0$ on $\Gamma$ by hypothesis.
Moreover, if our claim is not satisfied then, by applying Weak Harnack inequality repeatedly as in \eqref{lip1} (and Step 1 above) we obtain
$$w \leq f(\frac r2) \, \, e^{-c_0 r^{1-\beta}} \quad \text{in $\mathcal C_{r/2}(y',g(y'))$}.$$
 hence $$\frac 12 f(r) \le w(y) \le f(\frac r2) e^{-c_0 r^{1-\beta}},$$
and we reach a contradiction.

\qed

\begin{rem}\label{r1}
We remark that in Step 2 of Theorem \ref{main} when $\Gamma$ is Lipschitz we could have argued also as above, by constructing a sequence of points $y_k$ with $u(y_k) \to \infty$. This is precisely the strategy that appears in \cite{CFMS} for obtaining an upper bound on $u$.
\end{rem}

\subsection{NTA domains.} In \cite{JK}, the authors extended Theorem \ref{main} to a class of domains called ``NonTangentially Accessible" (NTA), in the case of operators in divergence form. The case of operators in non-divergence form is treated with probabilistic techniques in \cite{BB1}. We recall the definition of NTA domains which preserve the key properties of Lipschitz domains.

\begin{defn} A bounded domain $D \subset \R^n$ is NTA with constants $M$ and $r_0$ if
\begin{enumerate}
\item (Corkscrew condition.) $\forall x \in \p D$, $r<r_0$, $\exists y=y_r(x) \in D \cap B_{Mr}(x)$ such that 
$$ B_{r}(y) \subset D;$$
\item (Density estimate of the complement.) $\forall x\in \p D, r<r_0,$
$$|B_r(x) \setminus D| \geq M^{-1} |B_r(x)|.$$ 
\item (Harnack chain) If $\eps>0, x_1,x_2 \in D$, $dist(x_i,\p D)>\eps, |x_1-x_2| <k \eps,$ there exists a sequence of $Mk$ overlapping balls included in $D$, of radius $\eps/M$ such that, the first one is centered at $x_1$ and last one at $x_2$, and such that the centers of two consecutive balls are at most $\eps/(2M)$ apart.
\end{enumerate}
\end{defn}

Our strategy applies to NTA domains with very small modifications. We briefly mention how to define the sets $\mathcal C_r$ and $\mathcal A_r$ in this setting. 

Let $D \subset \R^n$ be NTA with constants $M, r_0.$ Say, $0 \in \p D.$ We define
$$\mathcal C_r:= D \cap B_r,$$
and 
$$\mathcal A_r := D_r \setminus \{d_\Gamma > \delta  \, r\},$$
where $d_\Gamma$ represents the distance function to $\Gamma:=\p D$,
$$d_\Gamma(x):= dist(x, \p D).$$
 Here constants depending on $n, \lambda, \Lambda, M, r_0$ are called universal.

The proof of Step 1 follows exactly as in the Lipschitz case. In the proof of \eqref{202}, the cylinder $Q_{2\delta}(x_0)$ is replaced by the ball $B_{2\delta}(x_0).$ Then weak Harnack inequality can be applied in view of (ii) the density property of NTA domains. After the iteration, the conclusion is that 
$$w>0 \quad \text{in the set $\{\delta |x| < d_{\Gamma}(x)\},$}$$
and this statement can be used in Step 3 to obtain the desired claim.

The proof of Step 2 can also be obtained as in the Lipschitz case, by observing that  $d_\Gamma^{-c}$ for $c>0$ small universal, is integrable in view of the corkscrew condition (via a covering argument.) Alternatively, we can also argue as in Remark \ref{r1}.

\subsection{Slit domains.} We discuss a version of Boundary Harnack in slit domains which is relevant in non-local problems via the Caffarelli-Silvestre extension \cite{CS}. A ``slit domain" is domain of the type $B_1 \setminus \mathcal P$ with $\mathcal P $ a closed subset of $\Gamma$, where $\Gamma$ is a Lipschitz graph as in Section 2 (or H\"older domain as in Subsection 3.1). The proof of the previous section leads to the following result. 

\begin{thm} Let $u, v\ge 0$ vanish continuously on $\mathcal P$ and assume that 
 \begin{equation}\label{305}
 \mathcal L u=\mathcal L v=0 \quad \text{in $B_1 \setminus \mathcal P$}.
 \end{equation}
 Suppose further that $u(\frac 1 2 e_n)=u(-\frac 1 2 e_n)=1.$ Then 
 $$C^{-1} \min\{v(\frac 1 2 e_n), v(-\frac 1 2 e_n)\} \leq \frac v u \leq C \max \{v(\frac 1 2 e_n), v(-\frac 1 2 e_n)\}$$
with $C=C(n,\lambda, \lambda, L)$ independent of $\mathcal P.$
\end{thm}
 
 Indeed, in this case the cylinders $\mathcal C_r$ and the sets $\mathcal A_r$ are defined as
 $$\mathcal C_r= \{|x'| <r, \quad |h_\Gamma|< r\}, \quad \mathcal A_r:=\{ |x'| <r, \quad \delta r< |h_\Gamma| < r\},$$
 and the proof remains identical.
 
 The conditions on $u$ and $v$ at the point $-e_n/2 $ can be removed from the statement of the theorem when there is a chain of overlapping balls of radius $\mu$ connecting the points $\pm e_n/2$, provided the constant $C$ depends on $\mu$ as well. 
 
Another situation when this happens, which as mentioned above appears in non-local problems, is when $\Gamma=\{x_n=0\}$ and $u$, $v$ are symmetric with respect to $\Gamma$. We have the following corollary.
 
 \begin{cor} Assume that $u,v \ge 0$ are even in the $x_n$ variable and vanish continuously on a closed subset $ \mathcal P \subset \{x_n=0\}$. Assume that $u,v$ satisfy \eqref{305} and are normalized such that $u=v=1$ at $e_n/2$. Then
 $$C^{-1} \le \frac vu \le C, \quad \mbox {in $B_{1/2} \setminus \mathcal P$},$$
 with $C=C(n,\lambda, \lambda)$. Moreover, $\frac u v$ is uniformly H\"older continuous in $B_{1/2} \setminus \mathcal P$.
 
 \end{cor}


\begin{thebibliography}{9999}
\bibitem[A]{A} Ancona A., {\it Principe de Harnack a la frontiere et theoreme de Fatou pour un operateur
elliptique dons un domaine lipschitzien}, Ann. Inst. Fourier 28 (1978) 169--213.
\bibitem[BB1]{BB1} Bass R.F. and Burdzy K., {\it A boundary Harnack principle in twisted H\"older domains,}
Ann. Math. 134 (1991) 253--276.
\bibitem[BB2]{BB2} Bass R.F. and Burdzy, K.,  {\it The boundary Harnack principle for non-divergence form elliptic operators},  J. London Math. Soc. (2) 50 (1994), no. 1, 157--169.
\bibitem[BBB]{BBB} Banuelos R., Bass R.F. and Burdzy K.,  {\it H\"older Domains and The Boundary Harnack Principle}, Duke Math. J., 64, 195--200 (1991).
\bibitem[CC]{CC} Caffarelli L. and  Cabre X., {\it Fully nonlinear elliptic equations}, American Mathematical Society Colloquium Publications, 43. American Mathematical Society, Providence, RI, 1995. vi+104 pp. 

\bibitem[CFMS]{CFMS} Caffarelli L., Fabes E., Mortola S. and Salsa S., {\it Boundary behavior of non-negative solutions of elliptic operators in divergence form,} Indiana Math. J.,30, 621--640 (1981).

\bibitem[CS]{CS} Caffarelli L. and Silvestre L.,  {\it An extension problem related to the fractional Laplacian}, Comm. Partial Differential Equations 32 (2007), no. 7-9, 1245--1260.
\bibitem[FGMS]{FGMS} Fabes E.B., Garofalo N., Marin-Malave S. and Salsa S., {\it Fatou theorems for some
nonlinear elliptic equations}, Rev. Mat. Iberoamericana, 4 (1988) 227--252.
\bibitem[D]{D}Dahlberg B., {\it On estimates of harmonic measure,} Arch. Rational Mech. Anal. 65 (1977), 272--288.
\bibitem[DG]{DG} De Giorgi E., {\it Sulla differenziabilita e l'analicita degli integrali multipli regolari,} Mem. Acca., Sci. torino, S. III, Parte I, (1957), 25--43.
\bibitem[F]{F} Ferrari F., {\it On boundary behavior of harmonic functions in H\"older domains}, Journal of Fourier Analysis and Applications, 1998, Volume 4, Issue 4-5, pp 447--461 (1988).

\bibitem[GT]{GT} Gilbarg D. and Trudinger N., {\it Elliptic partial differential equations of second order}, Second edition. Grundlehren der Mathematischen Wissenschaften [Fundamental Principles of Mathematical Sciences], 224. Springer-Verlag, Berlin, 1983. xiii+513 pp.
\bibitem[JK]{JK} Jerison, D.S. and Kenig, C.E., {\it  Boundary Behavior of Harmonic Functions in Non-tangentially Accessible Domains,} Adv. Math.,46, 80--147 (1982).
\bibitem[K]{K} Kemper J. T., {\it A boundary Harnack principle for Lipschitz domains and the principle of positive singularities,} Comm. Pure Appl. Math. 25 (1972), 247--255.
\bibitem[KS]{KS} Krylov N. V. and Safonov M. V., {\it An estimate for the probability of a diffusion process hitting a set of positive measure,} (Russian) Dokl. Akad. Nauk SSSR 245 (1979), no. 1, 18--20. 
\bibitem[W]{W}Wu J.-M. G., {\it Comparison of kernel functions, boundary Harnack principle, and relative
Fatou theorem on Lipschitz domains,} Ann. Inst. Fourier Grenoble 28 (1978)
147--167.

\end{thebibliography}
\end{document}